\documentclass[12pt]{article}

\newcounter{theor}
\newtheorem{theorem}{Theorem}[theor]
\newtheorem{lemma}{Lemma}[theor]
\newtheorem{corollary}{Corollary}[theor]

\newtheorem{proposition}{Proposition}[theor]
\newenvironment{proof}[1][Proof]{\begin{trivlist}
\item[\hskip \labelsep {\bfseries #1}]}{\end{trivlist}}
\newenvironment{definition}[1][Definition]{\begin{trivlist}
\item[\hskip \labelsep {\bfseries #1}]}{\end{trivlist}}

\hsize by -\marginparwidth \oddsidemargin -4mm \evensidemargin
\oddsidemargin

\usepackage{graphicx}
\newlength{\graphwidth}
\newlength{\graphheight}
\usepackage{graphicx}

\begin{document}
\centerline{\bf \large THE $\alpha$-INVARIANTS ON TORIC FANO
MANIFOLDS}

\bigskip
\centerline{Jian Song}
\bigskip

\centerline{Department of Mathematics} \centerline{Columbia
University, New York, NY 10027}
\bigskip
\bigskip

\centerline{\bf \S 1. Introduction}
\bigskip
\noindent The global holomorphic invariant $\alpha_{G}(M)$
introduced by Tian\cite{T1}, Tian and Yau\cite{TY} is closely
related to the existence of K\"ahler-Einstein metrics. In his
solution to the Calabi conjecture, Yau\cite{Y} proved the
existence of a K\"ahler-Einstein metric on compact K\"ahler
manifolds with nonpositive first Chern class. K\"ahler-Einstein
metrics do not always exist in the case when the first Chern class
is positive, for there exist known obstructions such as the Futaki
invariant. For a compact K\"ahler manifold $M$ with positive first
Chern class, Tian\cite{T1} proved that $M$ admits a
K\"ahler-Einstein metric if $\alpha_{G}(M)>\frac{n}{n+1}$, where
$n=\dim M$. In the case of compact complex surfaces, he proved
that any compact complex surface with positive first Chern class
admits a K\"ahler-Einstein metric except $CP^2\#1\overline{CP^2}$
and $CP^2\#2\overline{CP^2}$\cite{T3}.

\bigskip
\noindent There have been many nice results on the classification
of toric Fano manifolds. Mabuchi discovered that if a toric Fano
manifold is K\"ahler-Einstein then the barycenter of the
polyhedron defined by its anticanonical divisor is at the origin.
V. Batyrev and E. Selivanova \cite{BS} estimate the lower bound of
$\alpha$-invariant of symmetric toric Fano manifolds which is
sufficient to show the existence of K\"ahler-Einstein metric.

\bigskip
\noindent In this paper, we apply the Tian-Yau-Zelditch expansion
of the Bergman kernel on polarized K\"ahler metrics to approximate
plurisubharmonic functions and obtain a formula to calculate the
$\alpha_G$-invariants of all toric Fano manifolds precisely. This
gives a generalization of the result by V. Batyrev and E.
Selivanova\cite{BS} and also this formula confirms the earlier
result \cite{SJ} on the estimates of $\alpha$ invariants on
$CP^2\#1\overline{CP^2}$ and $CP^2\#2\overline{CP^2}$.

\bigskip
\noindent Our main theorems are
\begin{theorem}
If X is a toric Fano manifold of complex dimension n
then\\
$(a)$ $\alpha_G(X)=1$ if X is symmetric, otherwise\\
$(b)$ $\alpha_G(X)=\frac{\min_{0\neq v\in S}
\frac{|w_v|}{|v|}}{1+\min_{0\neq v\in S} \frac{|w_v|}{|v|}}\leq
\frac{1}{2}.$
\end{theorem}

\begin{corollary}
If X is a toric Fano manifold, then X is symmetric if and only if
$\alpha_G(X)=1$.
\end{corollary}

\begin{theorem}
If X is a toric Fano manifold, then $\{\alpha_{m,G}(X)\}_{m\geq
1}$ is stationary. More precisely, $\alpha_{m, G}(X)=\alpha_G(X)$
if $m\geq m_0$, where $m_0$ is the least positive integer such
that $m_0 v$ is an integral point and $v$ is the minimizer of
$\min_{0\neq v\in S} \frac{|w_v|}{|v|}.$
\end{theorem}

\bigskip
\bigskip
\noindent $Acknowledgements$. The author deeply thanks his
advisor, Professor D.H. Phong for his constant encouragement and
help. He also thanks Professor Jacob Sturm and Professor Zhiqin Lu
for their suggestion on this work. This paper is part of the
author's future Ph.D. thesis in Math Department of Columbia
University.

\bigskip
\bigskip
\centerline{\bf \S 2. Notations}

\setcounter{theor}{2}
\bigskip
\noindent Let us first recall the definition of an invariant
$\alpha_G(X)$ introduced by Tian. Let $X$ be an $n$-dimensional
compact complex manifold with positive first Chern class $c_1(X)$
and $G$ a compact subgroup of $Aut(X)$. Choose a $G$-invariant
K\"ahler metric $g=g_{i\overline{j}}$ on $X$ such that
$\omega_g=\frac{\sqrt{-1}}{2\pi}\sum g_{i\overline{j}}dz_i\wedge
d\overline{z}_j\in c_1(X).$

\begin{definition} Let $P_G(X,g)$ be the set of all $C^2$-smooth
$G$-invariant real-valued functions $\varphi$ such that $\sup_X
\varphi =0$ and
$\omega_g+\frac{\sqrt{-1}}{2\pi}\partial\overline{\partial}\varphi>0.$
The $\alpha_G(X)$ invariant is defined as superemum of all
$\lambda>0$ such that $$ \int_X e^{-\alpha\varphi}\omega^n \leq
C(\alpha)$$ for all $\varphi\in P_G(X,g)$, where $C(\alpha)$ is a
positive constant depending only on $\alpha, g$ and $X$.

\end{definition}

\noindent Let $N$ be a lattice of rank n, $M= Hom(N, Z)$ the dual
lattice. $M_R=M\otimes_Z R$, $N_R=N\otimes_Z R$. Let
$X=X_{\Sigma}$ be a smooth projective toric $n-$fold defined by a
complete fan $\Delta$ of regular cones $\Delta \subset M_R$ and
denote $\Delta(i)$ the $i$-dimensional cone of $\Delta$. We put
$T=C^*=$\{$(t_1, t_2,...,t_n)|$$t_i\in C^*$\}. For $a\in M$ and
$b\in N$, we define $<a, b>\in Z$, $\chi^a\in Hom_{alg gp}(T,C^*)$
by
\begin{eqnarray*}
&&<a, b>=\sum_{i=1}^n a_i b_i,\\
&&\chi^a((t_1,...,t_n))=t_1^{a_1}t_2^{a_2}...t_n^{a_n}.
\end{eqnarray*}

\bigskip
\noindent For each $\rho\in \Delta(1)$, let $b_{\rho}$ denote the
uniquie fundamental generator of $\rho$. We now consider the
divisor $K=-\sum_{\rho\in\Delta(1)}D(\rho)$ on $X=X_{\Delta}.$ The
following theorem is due to Demazure\cite{Dam}.
\begin{theorem}
K is a canonical divisor of $X_{\Delta}$, and the
following are equivalent:\\
$(a)$ $X_{\Delta}$ is a toric fano manifold.\\
$(b)$ $-K$ is ample.\\
$(c)$ $-K$ is very ample. \\
$(d)$ $\Sigma_{-K}=$\{ $a\in M_R | <a, b_{\rho}>\leq 1$ for all
$\rho \in \Delta(1)$\} is an n-dimensional compact convex
polyhedron whose vertices are exactly \{$a_\tau| \tau\in
\Delta(n)$\}, where each $a_{\tau}$ denotes the unique element of
$M$ such that $<a_{\tau}, b>=1$ for all fundamental generators $b$
of $\tau$.
\end{theorem}

\noindent The maximal torus $T\subset Aut(X)$ acting on $X$ has an
open dense orbit $U \subset X$, so the normalizer $N(T) \subset
Aut(X)$ of $T$ has a natural action on $U$. Let $W(X)=N(T)/T$ and
we identify the maximal torus $T\subset Aut(X)$ with an open dense
orbit $U$ in $X$ by choosing an arbitrary point$x_0\in U$, then we
have the following splitting short exact sequence

$$ 1\rightarrow T\rightarrow N(T)\rightarrow W(X)\rightarrow 1,$$
i.e., an embedding $W(X)\hookrightarrow N(T)$. Denote by
$K(T)=(S^1)^n$ the maximal compact subgroup in $T$. We choose $G$
to be the maximal compact subgroup in $N(T)$ generated by $W(X)$
and $K(T)$, so that we have the short exact sequence
$$ 1\rightarrow K(T) \rightarrow G \rightarrow W(X) \rightarrow 1
.$$

\begin{proposition}
Let X=$X_{\Delta}$ be a smooth projective toric n-fold defined by
a complete regular polyhedral fan $\Delta$. Then the group W(X) is
isomorphic to the finite group of all symmetries of $\Delta$,
i.e., W(X) is isomorphic to a subgroup of $GL(N)$$(\simeq
GL(n,Z))$ consisting of all elements $\gamma \subset GL(N)$ such
that $\gamma (\Delta)=\Delta$.
\end{proposition}

\noindent {\bf Remark:} $W(X)$ is as well isomorphic to a subgroup
of $GL(M)$$(\simeq GL(n,Z))$ consisting of all elements $\gamma
\subset GL(M)$ such that $\gamma(\Sigma)=\Sigma.$

\begin{definition}
A toric n-fold $X$ is symmetric if the trivial character is the
only $W(X)$-invariant algebraic character of T, i.e.
$N^{W(X)}=\{\chi\in N |$ $g\chi=\chi$ for all $g\in
W(X)\}$=$\{0\}$.
\end{definition}

\begin{definition}
\noindent Let $S=\{v\in \partial\Sigma|$ $gv=v$ for all $g\in
W(X)\}$ be the stable points of $W(X)$ on the boundary of
$\Sigma$. If $S\neq\{0\}$ then for any $0\neq v\in S$, we define
$w_v$ related with $v$ by $w_v=\partial \Sigma \cap \{-tv|$
$t\geq0\}.$
\end{definition}

\noindent {\bf Remark:} It's easy to see $X$ is symmetric if and
only if $S={0}.$

\bigskip

\bigskip
\bigskip
\centerline{\bf \S 3. Holomorphic approximation of PSH}
\bigskip
\setcounter{theor}{3} \noindent In this section, we will employ
the technique in \cite{T2,Z} to obtain the approximation of
plurisubharmonic functions by logarithms of holomorphic sections
of line bundles. The Tian-Yau-Zelditch asymptotic expansion of the
potential of the Bergman metric is given by the following
theorem\cite{Z}.
\bigskip
\begin{theorem} Let M be a compact
complex manifold of dimension n and let $(L,h)\rightarrow M$ be a
positive Hermitian holomorphic line bundle. Let g be the K\"ahler
metric on M corresponding to the K\"ahler form $\omega
_{g}=Ric(h)$. For each $m\in N$, $h$ induces a Hermitian metric
$h_{m}$ on $L^{m}.$ Let
$\{S_{0}^{m},S_{1}^{m},...,S_{d_{m-1}}^{m}\}$ be any orthonormal
basis of $H^{0}(M,L^{m})$, $d_{m}=\dim H^{0}(M,L^{m}),$ with
respect to the inner product:
$$
(S_{1},S_{2})_{h_{m}}=\int_{M}h_{m}(S_{1}(x),S_{2}(x))dV_{g},
$$
where $dV_{g}=\frac{1}{n!}\omega _{g}^{n}$ is the volume form of
$g$. Then there is a complete asymptotic expansion:

$$
\sum\limits_{i=0}^{d_{m}-1}||S_{i}^{m}(x)||_{h_{m}}^{2}\sim
a_{0}(x)m^{n}+a_{1}(x)m^{n-1}+a_{2}(x)m^{n-2}+...
$$
for some smooth coefficients $a_{j}(x)$\ with $a_{0}=1$. More
precisely, for any k:
$$
\left| \left|
\sum\limits_{i=0}^{d_{m}-1}||S_{i}^{m}(x)||_{h_{m}}^{2}-
\sum_{j<R}a_{j}(x)m^{n-j}\right| \right| _{C^{k}}\leq
C_{R,k}m^{n-R}
$$
where $C_{R,k}$\ depends on $R,k$\ and the manifold $M$.
\end{theorem}

\bigskip

\noindent Let
\begin{eqnarray*}
\;\tilde{\;\omega _{g}} &=&\omega _{g}+\sqrt{-1}\partial
\overline{\partial} \phi
>0\ \\
\widetilde{h} &=&he^{-\phi }
\end{eqnarray*}

\noindent Let $\widetilde{h}_{m}$ be the induced Hermitian metric
of $\widetilde{h}$ on $L^{m}$, $\{\widetilde{S}
_{0}^{m},\widetilde{S}_{1,...,}^{m}\widetilde{S}_{d_{m}-1}^{m}\}$
be any orthonormal basis of $H^{0}(M,L^{m})$, where $d_{m}=\dim
H^{0}(M,L^{m}),$ with respect to the inner product

$$
(S_{1},S_{2})_{\widetilde{h}_{m}}=\int_{M}\widetilde{h}%
_{m}(S_{1}(x),S_{2}(x))dV_{\widetilde{g}}\;.
$$
By Theorem 3.4, we have

$$\sum\limits_{i=0}^{d_{m}-1}||\widetilde{S}_{i}^{m}(x)||_{\widetilde{h}
_{m}}^{2} =\left( \sum\limits_{i=0}^{d_{m}-1}||\widetilde{S}
_{i}^{m}(x)||_{h_{m}}^{2}\right) e^{-m\phi }.$$ Thus
$$
\phi -\frac{1}{m}\log \left(
\sum\limits_{i=0}^{d_{m}-1}||\widetilde{S}
_{i}^{m}(x)||_{\widetilde{h}_{m}}^{2}\right)=-\frac{1}{m}\log
\left(
\sum\limits_{i=0}^{d_{m}-1}||\widetilde{S}_{i}^{m}(x)||_{\widetilde{h}
_{m}}^{2}\right)
$$
As $m\rightarrow +\infty$, we obtain for any positive integer $R$
\begin{eqnarray*}
&&\frac{1}{m}\log \left( \sum_{j < R}
\widetilde{a}_j(x)m^{n-j}\right) \\
&=&\frac{1}{m}\log m^{n}(\sum_{j <R}
\widetilde{a}_{j}(x)m^{-j}) \\
&=&\frac{n}{m}\log m+\frac{1}{m}\log (1+O(\frac{1}{m}))\rightarrow
0
\end{eqnarray*}

\noindent Thus we have the following corollary of the
Tian-Yau-Zelditch expansion.

\begin{corollary}
$$
\left\| \phi -\frac{1}{m}\log \left(
\sum\limits_{i=0}^{d_{m}-1}||\widetilde{
S}_{i}^{m}(x)||_{h_{m}}^{2}\right) \right\| _{C^{k}}\rightarrow
0,\;as\;m\rightarrow +\infty.
$$
\end{corollary}
In other words, any plurisubharmonic function can be approximated
by the logarithms of holomorphic sections of $L^{m}$.

\bigskip
\bigskip
\centerline{\bf \S 4. Proof of The Main Theorem }
\bigskip
\setcounter{theor}{4} \noindent Suppose $X_{\Delta}$ is Fano, then
one obtains a convex $W(X)$-invariant polyhedron $\Sigma$ in $M_R$
defined by $\Sigma$ = $\{$ $a \in M_R$ $|$ $<a, b_{\rho}> \leq 1$,
for all $\rho \in \Delta(1)$ $\}$ where $b_{\rho}$ is the
fundamental generator of $\rho$. Let $L(\Sigma)=\{v_0, v_1,...,
v_k\}=M\cap\Sigma$. Then $v_0, v_1,..., v_k$ determine algebraic
characters $\chi_i: T \rightarrow C^*$ of $T$(i=0, 1,..., k).
Moreover, we have
$$ |\chi_i(x)|^2=e^{<v_i, y>}, i=0, ..., k,$$where y is the image of
$x$ under the canonical projection $\pi : T \rightarrow M_R$. Let
us define $u: U \rightarrow R$ as follows:
$$ u=\log(\sum_{i=0}^k |\chi_i(x)|^2), x\subset U\simeq T.$$

\bigskip

\noindent Since $u$ is $K(T)$-invariant, $u$ descends to a
function $\tilde{u}: M_R \rightarrow R$ defined as
$$ \tilde{u}= \log(\sum_{i=0}^k e^{<v_i, y>}), y\subset M_R.$$
Consider the $G$-invariant hermitian metric
$g={g_{i\overline{j}}}$ on $X$ such that the restriction of the
corresponding to $g$ differential 2-form on $U$ is defined by
$$ \omega_g= \partial\overline{\partial} u.$$
The metric $g$ is exactly the pull-back of Fubini-Study metric
from $P^m$ with respect to the anticanonical embedding
$X\hookrightarrow P^m$ defined by the algebraic characters
$\chi_0, \chi_1, ..., \chi_k$.

\bigskip

\noindent Let $\Sigma^{(m)}$=$\{$ $a \in M_R$ $|$ $<a, b_{\rho}>
\leq m$ and $L(\Sigma^{(m)})=\{v_0, ...,v_{k_m}\}=M\cap
\Sigma^{(m)}$, where $k_m+1=\dim H^0(X, O((-K)^{m}))$ and
$\chi^{\mu}$: $T \rightarrow C^*$ defined by
$|\chi^{\mu}(x)|^2=e^{<\mu, y>}.$ We have the following lemma (see
\cite{F} p66)

\begin{lemma}
$H^0(X, O((-K)^m))=\oplus_{\mu\in L(\Sigma^{(m)})} C\cdot
\chi^{\mu}.$

\end{lemma}

\begin{proposition}
$\{\chi^{\mu}\}_{\mu\in L(\Sigma^{(m)})}$ is an orthogonal basis
of $H^0(X,O((-K)^(m))$ with respect to the inner product
$<,>_{h^m}$, where $h^m=\frac{1}{(\sum_{i=0}^k |\chi_i(x)|^2)^m}.$
\end{proposition}

\begin{proof}

\begin{eqnarray*}
&&\int_X <\chi^{\mu}, \chi^{\nu}>_{h^m}\omega^n\\
&=&\int_T \frac{(z_1^{\mu_1}...z_n^{\mu_n})(\overline{z}_1^{\nu_1}...\overline{z}_n^{\nu_n})}
{(\sum_{i=0}^k |z^{v_i}|^m)}\omega^n\\
&=&\int_T \frac{|z_1|^{\mu_1+\nu_1}...|z_n|^{\mu_n+\nu_n}e^{i(\mu_1-\nu_1)\theta_1}...e^{i(\mu_n-\nu_n)\theta_n}}
{(\sum_{i=0}^k |z^{v_i}|^m)}\omega^n.\\
\end{eqnarray*}
which is $0$ if $\mu\neq\nu$.

\end{proof}

\bigskip
\noindent For any $\varphi \in P_G(X, \omega)$, by Corollary 3.3
$$\varphi(x)= \lim_{m\rightarrow \infty} \frac{1}{m} \log
\frac{\sum_{\mu\in L(\Sigma^{(m)})} a_{\mu}^{(m)}
|\chi^{\mu}(x)|^2}{(\sum_{i=0}^k |\chi_i(x)|^2)^k}$$

\bigskip
\begin{lemma}
There exists $\epsilon > 0$ such that for $\varphi \in P_G(X,
\omega)$ and $\tilde{m}>0$ there exist $m>\tilde{m}$ and $\mu \in
L(\Sigma^{(m)})$ with $(a_{\mu}^{(m)})^{\frac{1}{m}}>\epsilon.$
\end{lemma}

\begin{proof}
Otherwise, for any $\epsilon>0$ there exists $\varphi_{\epsilon}$
and $\tilde{m}$ such that for any $m>\tilde{m}$ and $\mu\in
L(\Sigma^{(m)})$ we have $(a_{\mu}^{(m)})^{\frac{1}{m}}<\epsilon.$
By choosing $m$ large enough we have

\begin{eqnarray*}
\varphi_{\epsilon}(x) & \leq & \frac{1}{m} \log \frac{\sum_{\mu\in
L(\Sigma^{(m)})} |\chi^{\mu}(x)|^2}{(\sum_{i=0}^k
|\chi_i(x)|^2)^m}+\log\epsilon \\
&=&\frac{1}{m} \log (\sum_{\mu\in L(\Sigma^{(m)})}
\frac{|\chi^{\mu}(x)|^2}{(\sum_{i=0}^k
|\chi_i(x)|^2)^m})+\log\epsilon\\
&\leq& \frac{1}{m} \log (\sum_{\mu\in L(\Sigma^{(m)})} 1)+\log\epsilon\\
&\leq& Const+\log\epsilon.
\end{eqnarray*}

\noindent
Since $\epsilon$ can be chosen arbitrarily small, the
above inequality implies that $\varphi_{\epsilon} \rightarrow
-\infty$ uniformly as $\epsilon$ goes to $0$, which contradicts
the fact that $\sup_{X} \varphi= 0$.

\end{proof}

\bigskip
\noindent For any $\varphi\in P_G(X,\omega)$, by Lemma 4.1 we have
\begin{eqnarray*}
\varphi(x)&=&\lim_{m\rightarrow\infty} \frac{1}{m}\log
\frac{\sum_{\mu\in
L(\Sigma^{(m)})}a_{\mu}^{(m)}|\chi^{\mu}(x)|^2}{(\sum_{i=0}^k|\chi_i(x)|^2)^m}\\
&\geq& \frac{1}{m}\log \frac{\sum_{g\in W(X)}|\chi^{g\mu}(x)|^2}{(\sum_{i=0}^k|\chi_i(x)|^2)^m}-C_1\\
&\geq& \log \frac{|\chi^{\sum_{g\in W(X)}
g\mu}(x)|^{\frac{2}{m|W(X)|}}}{(\sum_{i=0}^k|\chi_i(x)|^2)}-C_1\\
&\geq& \log \frac{|\chi(x)|^{\frac{2\sum_{g\in W(X)}
g\mu}{m|W(X)|}}}{(\sum_{i=0}^k|\chi_i(x)|^2)}-C_1\\
\end{eqnarray*}

\noindent Put $v= \frac{\sum_{g\in W(X)}g\mu}{|W(X)|K}$, then we
have $\tilde{\varphi}(y) \geq \log\frac{e^{<v, y>}}{\sum_{i=0}^{m}
e^{<v_i, y>}}$

\noindent
Put $y_i=\log|t_i|^2$
$t_i=e^{\frac{yi}{2}+\sqrt{-1}\theta_i}$, then
\begin{eqnarray*}
&&\frac{dt_i}{t_i}=\frac{1}{2}dy_i+\sqrt{-1}d\theta_i\\
&&\frac{d\overline{t}_i}{\overline{t_i}}=\frac{1}{2}dy_i+\sqrt{-1}d\theta_i\\
&&\frac{dt_i\wedge d\overline{t}_i}{|t_i|^2}=-\sqrt{-1}dy_i\wedge
d\theta_i\\
&&\partial\overline{\partial}u=\sum_{i,j}\frac{\partial^2
u}{\partial y_i \partial y_j} \frac{dt_i\wedge
d\overline{t}_j}{t_i \overline{t}_j}\\
&&(\partial\overline{\partial}u)^n=\det(\frac{\partial^2u}{\partial
y_i\partial y_j})dy_1\wedge ...\wedge dy_n\wedge
d\theta_1\wedge...\wedge d\theta_n
\end{eqnarray*}

\begin{lemma}
Let
$\tilde{F}=e^{\tilde{u}}\det\frac{\partial^2\tilde{u}}{\partial
y_i\partial y_j}$, then $0<c\leq \tilde{F}\leq C$.
\end{lemma}

\begin{proof}

\noindent $e^{-u}\frac{dt_1\wedge d\overline{t}_1\wedge...\wedge
dt_n\wedge
d\overline{t}_n}{|t_1|^2...|t_n|^2}=e^{-\tilde{u}}dy_1\wedge...\wedge
dy_n\wedge d\theta_1\wedge...\wedge d\theta_n$
 can be extended to a
non-vanishing volume form on $X$. Also
\begin{eqnarray*}
(\partial\overline{\partial}u)^n&=&\det\frac{\partial^2 u
}{\partial t_i\overline{
\partial} t_j}dt_1\wedge d\overline{t}_1\wedge...\wedge
dt_n\wedge
d\overline{t}_n\\
&=&\det\frac{\partial^2\tilde{u}}{\partial y_i\partial
y_j}dy_1\wedge...\wedge dy_n\wedge d\theta_1\wedge...\wedge
d\theta_n
\end{eqnarray*}
is a non-vanishing volume form, so the quotient of these two
volume form must be positive and bounded. which proves the lemma.
\end{proof}

\bigskip\noindent
Now we can prove the Theorem 1.1.
\begin{eqnarray*}
\int_X e^{-\alpha\varphi}\omega^n&=&\int_X
e^{-\alpha\varphi}(\partial\overline{\partial}u)^n\\
&=&\int_{R^n}
e^{-\alpha\tilde{\varphi}}\det\frac{\partial^2\tilde{u}}{
\partial y_i \partial y_j} dy_1...dy_n\\
&\leq&\int_{R^n} e^{-\alpha\tilde{\varphi}-\tilde{u}}dy_1...dy_n\\
&\leq&\int_{R^n} e^{-\alpha \log \frac{e^{<v,y>}}{\sum
e^{<v_i,y>}}-\log(\sum e^{<v_i,y>})}dy_1...dy_n\\
&=&\int_{R^n} \frac{e^{-\alpha <v,y>}}{(\sum e^{<v_i,
y>})^{1-\alpha}}dy_1...dy_n
\end{eqnarray*}

\bigskip

\bigskip
\noindent If the stable points $S=\{0\}$, then $X$ is symmetric so
that $v= \frac{\sum_{g\in W(X)}g\mu}{m|W(X)|}=0$ for all $\mu\in
L(\Sigma^{(m)})$. Therefore for all $\alpha<1$ the integral
$$\int_X e^{-\alpha\varphi} \omega^n=\int_{R^n} \frac{1}{(\sum_{i=0}^k e^{<v_i,
y>})^{1-\alpha}}dy_1...dy_n$$ is finite since every
$n$-dimensional cone $\sigma_j\in \Delta$$(j=1,...,l)$ is
generated by a basis of the lattice $N$ and the fact that
$N_R=\sigma_1\cup...\cup \sigma_l.$ This implies $\alpha_G(X)\geq
1$ so that by Tian's theorem\cite{T1} $X$ admits Kahler-Einstein
metric. This is a result by V. Batyrev and E.N.
Selivanova\cite{BS}.

\bigskip
\noindent If $S\neq\{0\}$ then for any $0\neq v\in S$, we have
$w_v\in\partial\Sigma$ related with $v$ satisfying $$<w_v,
v>=-|w_v||v|.$$

\noindent The integral $$\int_{R^n} \frac{e^{-\alpha <v,y>}}{(\sum
e^{<v_i, y>})^{1-\alpha}}dy_1...dy_n=\int_{R^n}
\left(\frac{e^{-\frac{\alpha}{1-\alpha} <v,y>}}{\sum e^{<v_i,
y>}}\right)^{1-\alpha}dy_1...dy_n$$ is finite if
$$-\frac{\alpha}{1-\alpha}v\in int({\Sigma}).$$
i.e.
$$<-\frac{\alpha}{1-\alpha}v,w_v>=
\frac{\alpha}{1-\alpha}|v||w_v|\leq|w_v|^2$$

\bigskip

\noindent
 Then for all $\alpha < \frac{\min_{0\neq v\in S}
\frac{|w_v|}{|v|}}{1+\min_{0\neq v\in S} \frac{|w_v|}{|v|}}$ the
integral $\int_X e^{-\alpha\varphi} \omega^n$ is finite. Therefore
$$\alpha_G(X)\geq\frac{\min_{0\neq v\in S}
\frac{|w_v|}{|v|}}{1+\min_{0\neq v\in S} \frac{|w_v|}{|v|}}.$$

\bigskip

\noindent In order to estimate the upper bound of $\alpha_G(X)$ we
will construct a sequence of PSH functions. Suppose $S\neq\{0\}$,
then for all $\alpha$ with $1>\alpha > \frac{\min_{0\neq v\in S}
\frac{|w_v|}{|v|}}{1+\min_{0\neq v\in S} \frac{|w_v|}{|v|}} $ we
choose
$\tilde{\varphi_{\epsilon}}=\log(\frac{e^{<\tilde{v},y>}+\epsilon}{\sum
e^{<v_i, y>}})$ which is increasing and uniformly bounded from
above where $\min_{0\neq v\in S} \frac{|w_v|}{|v|}$ is achieved at
$\tilde{v}\in S$ . Then by Fatou lemma we have
$$\lim_{\epsilon\rightarrow 0}\int_X e^{-\alpha
\varphi_{\epsilon}}\omega^n=\int_X
e^{-\alpha\varphi_{0}}\omega^n=\infty$$ which implies
$\alpha_G(X)\leq \frac{\min_{0\neq v\in S}
\frac{|w_v|}{|v|}}{1+\min_{0\neq v\in S} \frac{|w_v|}{|v|}}$.
Combined the above estimates together, we have proved the Theorem
1.1.

\bigskip
 \noindent Also it's obvious to see that $\min_{0\neq v\in
S} \frac{|w_v|}{|v|}\leq 1$ for non-symmetric toric Fano manifold
$X$ thus $\alpha_G(X)\leq \frac{1}{2}$. This shows that there
doesn't exist any non-symmetric toric Fano manifold such that its
$\alpha_G$-invariant is greater than $\frac{n}{n+1}$ which is a
sufficient condition for the existence of K\"ahler-Einstein
metrics.

\bigskip
\noindent Now we prove Theorem 1.2 which is a direct corollary of
the proof of Theorem 1.1. \noindent Define
$P_{m,G}(X)$=$\{\varphi\in C^{\infty}(X,R)| \sup_X \varphi=0,$
$\varphi$ is $G$-invariant and there exists a basis $\{
S_i^{m}\}_{1\leq i\leq N_m}$ of $H^0(X, K_X^{-m})$ such that
$\omega_g
+\partial\overline{\partial}\varphi=\frac{1}{m}\partial\overline{\partial}\log(\sum_{i=0}^{N_m}|S_i^m|^2)\},$
where $N_m+1=\dim H^0(X,K_X^{-m})$ and $m$ is large.

\bigskip
 \noindent We also define for $m$ large,
$\alpha_{m,G}(X)=\sup\{\alpha |$  there exists $C>0$ such that for
all $\varphi\in P_{m,G}(X)$, $\int_X e^{-\alpha \varphi}dV\leq
\infty \}$.
\bigskip

\noindent It's easy to see $\alpha_{m,G}(X)$ is decreasing as $m$
goes to the infinity. By the argument to give the upper bound for
the $\alpha_G$-invariant, we can directly have the following
corollary which answers the question proposed by Tian\cite{T2} in
the special case of toric Fano manifolds.

\begin{corollary}
If X is a toric Fano manifold, then $\{\alpha_{m,G}(X)\}_m$ is
decreasing and stationary. More precisely, $\alpha_{m,
G}(X)=\alpha_G(X)$ if $m\geq m_0$, where $m_0$ is the least
positive integer such that $m_0 v\in M$ and $v$ is the minimizer
of $\min_{0\neq v\in S} \frac{|w_v|}{|v|}.$
\end{corollary}

\bigskip
\centerline{\bf \S 5. Multiplier ideal sheaf }
\setcounter{theor}{5}
\bigskip
\noindent In this section we relate the $\alpha$-invariant on
toric Fano manifolds with the method of the multiplier ideal sheaf
employed by Nadel\cite{N}. Here we follow the lines in \cite{DK}.
\begin{theorem}\textnormal{(Nadel)}
Let X be a Fano manifold of dimension n and G be a compact
subgroup of the group of complex automorphisms of X. Then X admits
a G-invariant K\"ahler-Einstein metric, unless $K^{-1}_X$
possesses a G-invariant singular hermitian metric
$h=h_0e^{-\varphi}$($h_0$ is a smooth G-invariant metric and
$\varphi$ is a G-invariant function in $L^1_{loc}(X))$, such that
the following properties occur.
\begin{enumerate}
\item h has a semipositive curvature current
   $$\Theta_h= -\frac{i}{2\pi}\partial\overline{\partial} \log
   h=\Theta_{h_0}+
   \frac{i}{2\pi}\partial\overline{\partial}\varphi\leq 0.$$
\item For every $\gamma\in(\frac{n}{n+1},1)$, the multiplier ideal
sheaf $\mathcal{J}(\gamma\varphi)$ is nontrivial, (i.e. $0\neq
\mathcal{J}(\gamma\varphi)\neq O_X$).
\end{enumerate}
\end{theorem}

\begin{theorem}\textnormal{(Nadel)}
Let $(X,\omega)$ be a K\"ahler manifold and let L be a holomorphic
line bundle over X equipped with a singular hermitian metric h of
weight $\phi$ with respect to a smooth metric $h_0$(i.e.
$h=h_0e^{-\phi})$. Assume that the curvature form $\Theta_h(L)$ is
positive definite in the sense of currents, i.e.
$\Theta_h(L)\geq\epsilon\omega$ for some $\epsilon>0$. Then we
have $H^q(X,K_X \otimes L\otimes \mathcal{J}(\phi))=0$   for all
$q\geq 1$
\end{theorem}

\begin{corollary}\textnormal{(Nadel)}
Let X, G, h and $\varphi$ be in theorm 4.1. then for all
$\gamma\in (\frac{n}{n+1}, 1)$,
\begin{enumerate}
\item The multiplier ideal sheaf $\mathcal{J}(\gamma\varphi)$ satisfies
$H^q(X, \mathcal{J}(\gamma\varphi)=0$ for all $q\geq1$.

\item The associated subscheme $V_{\gamma}$ of structure sheaf
$O_{V_{\gamma}}=O_X/\mathcal{J}(\gamma\varphi)$ is nonempty,
distinct from X, G-invariant and satisfies $H^q(V_{\gamma},
O_{V_{\gamma}})$= C for q=0 and vanishes for $q\geq 1$.
\end{enumerate}
\end{corollary}

\noindent In order to construct Kahler-Einstein metrics it's
sufficient to rule out the existence of any $G$-invariant
subscheme with the property (2) in the corollary.

\noindent In the case of toric Fano manifolds, we have the
following theorem.

\begin{theorem}
Let $X$ be a toric Fano manifold. If $X$ is not symmetric then
there always exists a $G$-invariant subscheme with the property
(2) in the corollary.
\end{theorem}

\begin{proof}
If $X$ is not symmetric then $\alpha_G(X)\leq \frac{1}{2}$ and we
can construct $G$-invariant $\varphi\in L^1_{loc}(X)$ such that
for all $\gamma\in (\alpha_G(X), 1)$
$$ \int_X e^{-\gamma\varphi}\omega^n=+\infty$$
therefore $J(\gamma\varphi)$ is nontrivial and there exist
subschemes $V_{\gamma}$ satisfying property $(2)$ of the
corollary.
\end{proof}

\bigskip

\centerline{\bf \S 6. Examples}
\bigskip
\noindent In this sections we will calculate the $\alpha$
invariants for $2$-dimensional toric Fano manifolds.

\bigskip
 \noindent
 Here $(1)$ $(2)$ $(3)$ $(4)$ are corresponding to $CP^2$
and $CP^2$ blow-up at $1$, $2$ and $3$ points and $(5)$ $(6)$
$(7)$ $(8)$ are the corresponding polyhedrons.

\noindent $CP^2$ and $CP^2$ blow-up at 3 points are symmetric thus
the their $\alpha_G$-invariants are both equal to 1.

\bigskip
\noindent For $CP^2\#1\overline{CP^2}$ its stable points of $G$ on
the boundary of the polyhedron in $(6)$ is $(\frac{1}{2},
\frac{1}{2})$ and $(-\frac{1}{2},-\frac{1}{2})$ and
$\frac{|(1/2,1/2)|}{|(-1/2,-1/2)|}=1$ then it is easy to see
$\alpha_G(CP^2\#1\overline{CP^2})=\frac{1}{2}.$

\bigskip
\noindent For $CP^2\#2\overline{CP^2}$ its stable points of $G$ on
the boundary of the polyhedron in $(7)$ is $(\frac{1}{2},
\frac{1}{2})$ and $(-1,-1)$ and
$\frac{|(1/2,1/2)|}{|(-1,-1)|}=\frac{1}{2}$ then it is easy to see
$\alpha_G(CP^2\#2\overline{CP^2})=\frac{1}{3}.$

\bigskip
\noindent The above calculation confirms our earlier
results\cite{SJ}.

\clearpage \setlength{\graphwidth}{450bp}
\setlength{\graphheight}{130bp}

\begin{figure}[p]
\begin{center} \footnotesize

  \includegraphics[width=\graphwidth,height=\graphheight]
     {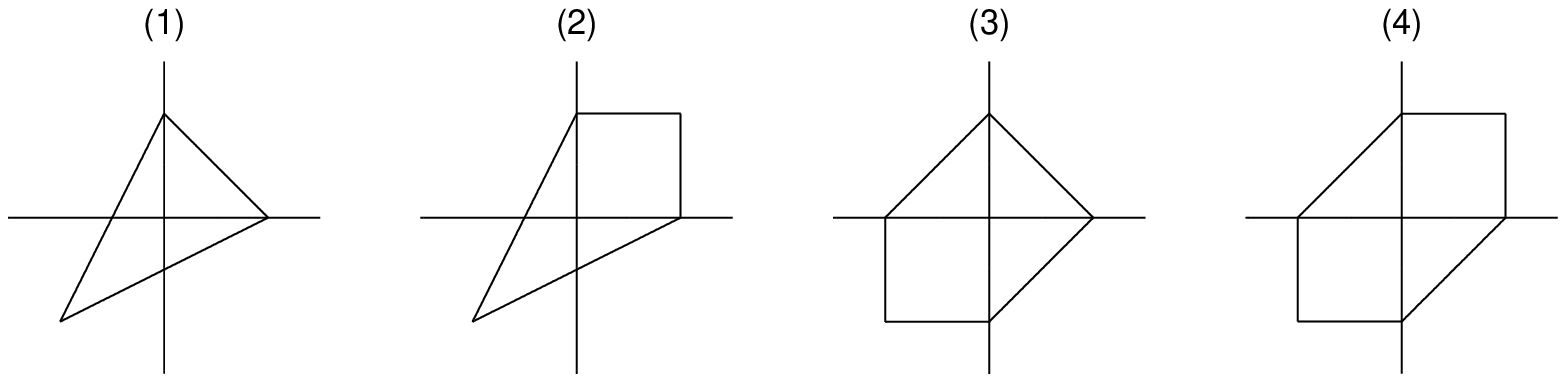}
\includegraphics[width=\graphwidth,height=\graphheight]
     {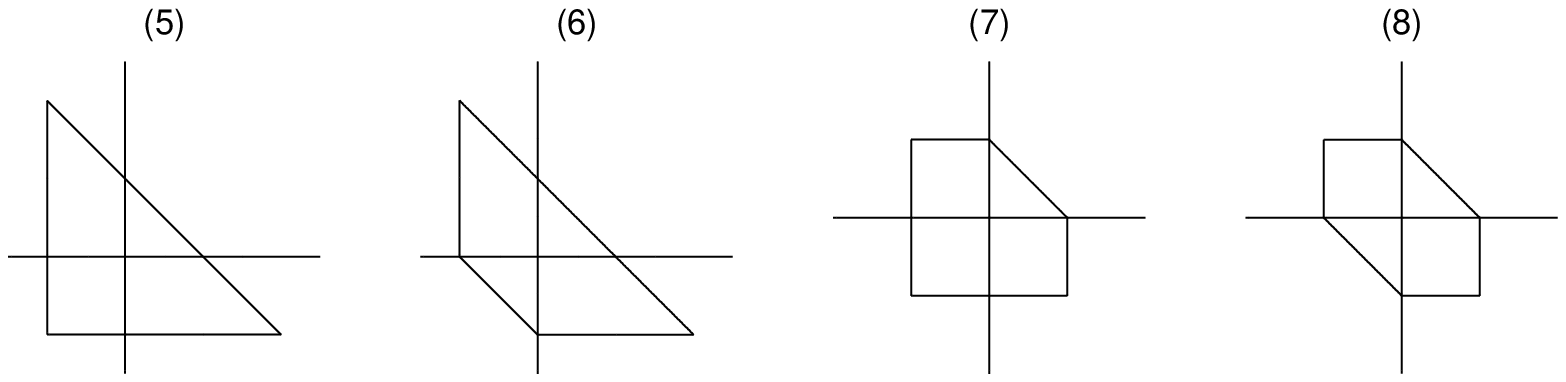}
\end{center}
\begin{center}

\end{center}
\end{figure}

\end{document}